\documentclass[12pt]{amsart} 

\title
[Boundary Coinvariants for Affine Buildings]
{Torsion in Boundary Coinvariants and K-theory for Affine Buildings}
   
\author{Guyan Robertson}        
\date{January 15, 2005}          
\address{School of Mathematics and Statistics, University of Newcastle, NE1 7RU, U.K.}
\email{a.g.robertson@newcastle.ac.uk}
 
\thanks{ \hfill Typeset by  \AmS-\LaTeX}

\usepackage{amsmath} 
\usepackage{amscd}

\input{prepictex} 
\input{pictex} 
\input{postpictex} 
\chardef\bslash=`\\ 
 
\makeatletter 
\def\verbatim{\interlinepenalty\@M \@verbatim 
  \leftskip\@totalleftmargin\advance\leftskip2pc 
  \frenchspacing\@vobeyspaces \@xverbatim} 
\makeatother 
\newtheorem{theorem}{Theorem}[section] 
\newtheorem{corollary}[theorem]{Corollary} 
\newtheorem{lemma}[theorem]{Lemma} 
\newtheorem{proposition}[theorem]{Proposition}

\theoremstyle{definition} 
\newtheorem{definition}[theorem]{Definition} 
\newtheorem{remark}[theorem]{Remark} 
 
\newtheorem{example}[theorem]{Example}

\newcounter{picture} 

\newsymbol\rtimes 226F 

\DeclareMathOperator{\covol}{covol}
\DeclareMathOperator{\ord}{ord}
\DeclareMathOperator{\rank}{rank}

\newcommand{\CC}{{\mathbb C}} 
\newcommand{\FF}{{\mathbb F}} 
 
\newcommand{\KK}{{\mathbb K}}

\newcommand{\QQ}{{\mathbb Q}}

\newcommand{\ZZ}{{\mathbb Z}} 
\newcommand{\cA}{{\mathcal A}}

\newcommand{\cS}{{\mathcal S}}

\newcommand{\D}{{\Delta}} 
 
\newcommand{\G}{{\Gamma}}

\newcommand{\fa}{{\mathfrak a}}
\newcommand{\fb}{{\mathfrak b}} 

\newcommand{\fD}{{\mathfrak D}}
 
\newcommand{\fs}{{\mathfrak s}} 
\newcommand{\ft}{{\mathfrak t}}

\newcommand{\SL}{{\text{\rm{SL}}}}

\newcommand{\pgl}{{\text{\rm{PGL}}}}
\newcommand{\psl}{{\text{\rm{PSL}}}}

\newcommand{\1}{{\bf 1}}

\begin{document} 

\begin{abstract}
Let $(G,{\mathfrak I},N,S)$ be an affine topological Tits system, and let 
$\Gamma$ be a torsion free cocompact lattice in $G$. This article studies the coinvariants
$H_0(\Gamma; C(\Omega,{\mathbb Z}))$, where $\Omega$ is the Furstenberg boundary of $G$.
It is shown that the class $[\1]$ of the identity function in $H_0(\Gamma; C(\Omega,{\mathbb Z}))$  has 
finite order, with explicit bounds for the order. 

A similar statement applies to the $K_0$ group of the boundary crossed product $C^*$-algebra $C(\Omega)\rtimes\Gamma$. If the Tits system has type $\widetilde A_2$, exact computations are given, both for the crossed product algebra and for the reduced group $C^*$-algebra.
\end{abstract}

\maketitle

\section{Introduction}

This article is concerned with coinvariants for group actions on the boundary of
an affine building. The results are most easily stated for subgroups of linear algebraic groups.
Let $k$ be a non-archimedean local field with finite residue field $\overline k$ of order $q$. Let $G$ be the group of $k$-rational points of an absolutely almost simple, simply connected linear algebraic $k$-group. Then $G$ acts on its Bruhat-Tits building $\Delta$, and on its Furstenberg boundary $\Omega$. 

Let $\Gamma$ be a torsion free lattice in $G$. The abelian group $C(\Omega,\ZZ)$ of continuous integer-valued functions on $\Omega$ has the structure of a $\Gamma$-module.  The module of $\Gamma$-coinvariants $\Omega_{\Gamma}= H_0(\Gamma; C(\Omega,\ZZ))$ is a finitely
generated group. We prove that the class $[\1]$  in  $\Omega_{\Gamma}$ of the constant function $\1\in C(\Omega,\ZZ)$ has finite order.
If $G$ is not one of the exceptional types $\widetilde E_8, \widetilde F_4$ or $\widetilde G_2$, 
then the order of $[\1]$ is less than $\covol(\Gamma)$, where the Haar measure $\mu$ on $G$ is normalized so that an Iwahori subgroup of $G$ has measure 1.
There is a weaker estimate for groups of exceptional type. If $G$ has rank 2 then the estimates are significantly improved.

The topological action of $\Gamma$ on the Furstenberg boundary is encoded in the crossed product $C^*$-algebra $\cA_\Gamma=C(\Omega)\rtimes\Gamma$.
Embedded in $\cA_\Gamma$ is the reduced group $C^*$-algebra $C^*_r(\G)$, which is the completion  of the complex group algebra of $\G$ in the regular representation as operators on $\ell^2(\G)$.
The action of $\Gamma$ on $\Omega$ is amenable, so the K-theory of $\cA_\Gamma$
is computable by known results, in contrast to that of $C^*_r(\G)$, which rests on the validity of the Baum-Connes conjecture.
The natural embedding  $C(\Omega) \to \cA_\Gamma$ induces a homomorphism 
$$\varphi:\Omega_\Gamma\to K_0(\cA_\Gamma)$$
and $\varphi([\1])=[\1]_{K_0}$, the class of \1 in the $K_0$-group of $\cA_\Gamma$.
Therefore $[\1]_{K_0}$ has finite order in $K_0(\cA_\Gamma)$.

If $\Gamma$ is a torsion free lattice in $G=\SL_3(k)$ then exact computations can be performed. The Baum-Connes Theorem of V. Lafforgue \cite{la} is used to compute $K_*(C_r^*(\G))$ and the results of \cite{RS} are used to compute $K_*(\cA_\Gamma)$.
In particular $K_0(C_r^*(\G)) = \ZZ^{\chi(\G)}$, a free abelian group,
one of whose generators is the class $[\1]$.
The embedding of $C_r^*(\G)$ into $\cA_\G$ induces a homomorphism $\psi :  K_*(C_r^*(\G)) \to K_*(\cA_\Gamma)$.  This homomorphism is not injective, since $[\1]$ has finite order in $K_0(\cA_\Gamma)$.
The computations at the end of the article suggest that the {\it only} reason for failure of injectivity of the homomorphism $\psi$ is the fact that $[\1]$ has finite order in $K_0(\cA_\Gamma)$.

Much of this article considers the more general case where $\Gamma$ is a subgroup of a topological group $G$ with a BN-pair, and $\Gamma$ acts on the boundary $\Omega$ of the affine building of $G$. 

The results are organized as follows. Sections \ref{one} and \ref{two} state and prove the main result concerning the class $[\1]$ in $\Omega_\Gamma$. 
Section \ref{rank2} gives improved estimates in the rank 2 case.
Section \ref{boundaryalgebra} studies the connection with the  K-Theory of the boundary algebra $\cA_\Gamma$. Comparison with K-theory of the reduced $C^*$-algebra $C^*_r(\Gamma)$ is made in Section \ref{a2tilde}, which contains some exact computational results for buildings of type
$\widetilde A_2$.

\section{Torsion in Boundary Coinvariants}\label{one}

Let $(G,{\mathfrak I},N,S)$ be an affine topological Tits system \cite[Definition 2.3]{garland}. Then  $G$ is a group with a BN-pair
in the usual algebraic sense \cite[Section 2]{ti1} and the Weyl group $W=N/({\mathfrak I}\cap N)$ is an infinite Coxeter group with generating set $S$.  The subgroup ${\mathfrak I}$ of $G$ is called an {\it Iwahori} subgroup. A subgroup of $G$ is {\it parahoric} if it contains a conjugate of ${\mathfrak I}$. The topological requirements are that $G$ is a second countable locally compact group and that all proper parahoric subgroups of $G$ are open and compact \cite[Definition 2.3]{garland}.

Let $n+1=|S|$ be the {\it rank} of the Tits system.
The group $G$ acts on the Tits complex $\Delta$, which is an affine building of dimension $n$.
It will be assumed throughout that $\Delta$ is irreducible; in other words, the Coxeter group $W$
is not a direct product of nontrivial Coxeter groups. Denote by $\Delta^i$ the set of $i$-simplices of $\Delta$, ($0\leq i\leq n$). The vertices of $\Delta$ are the maximal proper parahoric subgroups of $G$, and a finite set of such subgroups spans a simplex in $\Delta$ if and only if its intersection is parahoric.
The action of $G$ on $\Delta$ is by conjugation of subgroups. The building $\Delta$ is a union of $n$-dimensional subcomplexes, called {\it apartments}.
Each apartment is a Coxeter complex, with Coxeter group $W$.

Associated with the Coxeter system $(W,S)$ there is a Coxeter diagram of type $\widetilde X_n$  ($X=A,B,\dots ,G$), whose vertex set $I$ is a set of $n+1$ {\em types}, which are in natural bijective correspondence with the elements
of $S$. Each vertex $v\in\Delta^0$ has a {\em type} $\tau(v)\in I$. The type of a simplex in $\Delta$ is the set of types of its vertices. 
By construction, the action of $G$ on $\Delta$ preserves types.
A type $\ft \in I$ is {\it special} if deleting $\ft$ and all the edges containing $\ft$ from the diagram of type $\widetilde X_n$ results in the diagram of type $X_n$ (the diagram of the corresponding finite Coxeter group). A vertex $v\in \Delta$ is said to be {\em special} if its type $\tau(v)$ is special \cite[1.3.7 ]{bt}.

A simplex of maximal dimension $n$ in $\Delta$ is called a {\em chamber}. Every chamber has exactly one vertex of each type. 
If $\sigma$ is any chamber containing the vertex $v$ then the codimension-1 face of $\sigma$ which does not contain $v$ has type $I-\{\tau(v)\}$. 

The action of $G$ on $\Delta$ is {\it strongly transitive}, in the sense that $G$ acts transitively on the set of pairs $(\sigma,A)$ where $\sigma$ is a chamber contained in an apartment $A$ of $\Delta$. 
The building $\Delta$ is {\it locally finite}, in the sense that the number of chambers containing any simplex is finite, and it is {\em thick}, in the sense that each simplex of dimension $n-1$  is contained in at least three chambers. If $\tau$ is a simplex in $\Delta$ of dimension $n-1$ and type $I-\{\ft\}$, then the number of chambers of $\Delta$ which contain $\tau$ is $q_\ft+1$ where $q_\ft\ge2$. The integer $q_\ft$ depends only on $\ft$; not on $\tau$.

Associated with the group $G$ there is also a spherical building, the building at infinity $\Delta_\infty$.
The {\it boundary} $\Omega$ of $\Delta$ is the set of chambers of $\Delta_\infty$, endowed with a 
natural compact totally disconnected topology, which we shall describe later on. 
Since $G$ acts transitively on the chambers of  $\Delta_\infty$, $\Omega$ may be identified with the topological homogeneous space $G/B$, where the 
{\it Borel subgroup} $B$ is the stabilizer of a chamber of $\Delta_\infty$.

\begin{example}
A standard example is $G=\SL_{n+1}(\QQ_p)$, where $\QQ_p$ is the field of $p$-adic numbers. In this case $B$ is the subgroup of upper triangular matrices in $G$, and $\Omega$ is the Furstenberg boundary of $G$.
\end{example}

If $\Gamma$ is a subgroup of $G$, then $\Gamma$ acts on $\Omega$, and the abelian group $C(\Omega,\ZZ)$ of continuous integer-valued functions on $\Omega$ has the structure of a $\Gamma$-module.  The module of $\Gamma$-coinvariants, $C(\Omega,\ZZ)_\Gamma$, is the quotient of 
$C(\Omega,\ZZ)$ by the submodule generated by 
$\{g\cdot f-f : g\in\Gamma, f\in C(\Omega,\ZZ)\}$. 
Recall that $C(\Omega,\ZZ)_\Gamma$ is the homology group  $H_0(\Gamma; C(\Omega,\ZZ))$. For the rest of this article, $C(\Omega,\ZZ)_\Gamma$ will be denoted
simply by $\Omega_\Gamma$.
Define $c(\Gamma)\in \ZZ_+\cup\{\infty\}$ to be the number of $\Gamma$-orbits of chambers in $\Delta$.

If $\Gamma$ is a torsion free cocompact lattice in $G$, then  $c(\Gamma)$ is the number of $n$-cells of the finite
cell complex $\Delta\backslash\Gamma$. 
Suppose that the Haar measure $\mu$ on $G$ has the Tits normalization $\mu({\mathfrak I})=1$ \cite[\S 3.7]{tit}.
Then $c(\Gamma)=\covol(\Gamma)$.

We shall see below that if $\Gamma$ is a torsion free
cocompact lattice in $G$ then  $\Omega_\Gamma$ is a finitely generated abelian group.
Note that such a torsion free lattice $\Gamma$ acts freely and properly on $\Delta$
\cite[Lemma 2.6, Lemma 3.3]{garland}.
If $f\in C(\Omega,\ZZ)$ then $[f]$ will denote its class in $\Omega_\Gamma$. 
Also, $\1$ will denote the constant function defined by $\1(\omega)=\omega$ for $\omega\in\Omega$. 

\begin{theorem}\label{main} Let $(G,{\mathfrak I},N,S)$ be an affine topological Tits system and
let $\Gamma$ be a torsion free lattice in $G$.  Then $\Omega_\Gamma$
is a finitely generated abelian group and the following statements hold.
\begin{itemize}
\item[(1)] The element $[\1]$ has finite order in $\Omega_\Gamma$.
\item[(2)] If $\fs\in I$ is a special type, then the order of $[\1]$ in $\Omega_\Gamma$ satisfies
\[
\ord([\1]) < q_\fs \cdot \covol(\Gamma) \, .
\]
\item[(3)]
If, in addition, $G$ is not one of the exceptional types $\widetilde G_2, \widetilde F_4, \widetilde E_8$, then
\[
\ord([\1]) < \covol(\Gamma) \, .
\]
\end{itemize}
\end{theorem}

\begin{remark}
A torsion free lattice in $G$ is automatically cocompact \cite[II.1.5]{ser}.
\end{remark}

\begin{remark}\label{supergroup}
Suppose that $\Gamma$ is isomorphic to a subgroup of a group $\Gamma'$ and that the action of $\Gamma$ on
$\Omega$ extends to an action of $\Gamma'$ on $\Omega$. Then there is a natural
surjection  $\Omega_\Gamma \to \Omega_{\Gamma'}$. It follows that Theorem \ref{main}
remains true if $\Gamma$ is replaced by any such group $\Gamma'$.
\end{remark}

\begin{remark} The group $\Omega_\Gamma$ depends only on $\Gamma$ and not on the ambient group $G$. 
This follows from the rigidity results of \cite{kl}, if $n\ge 2$, and from  \cite {gr}
if $n=1$.
\end{remark}

We now describe briefly how Theorem \ref{main} applies to algebraic groups.
Let $k$ be a non-archimedean local field and 
let $G$ be the group of $k$-rational points of an absolutely almost simple, simply connected linear algebraic $k$-group : e.g. $k=\QQ_p$,\, $G=\SL_{n+1}(\QQ_p)$. Associated with $G$ there is a topological Tits system  of rank $n+1$, where $G$ has $k$-rank $n$ \cite{IM}.
Now $G$ acts properly on the corresponding Bruhat-Tits building $\Delta$ \cite[\S 2.1]{tit}, and on the boundary $\Omega=G/B$, where $B$ is a Borel subgroup \cite[Section 5]{bm}.

Let $q$ be the order of the residue field $\overline k$. 
For each type $\ft\in I$ there is an integer $d(\ft)$ such that $q_\ft=q^{d(\ft)}$. That is, any simplex $\tau$ of codimension one and type $I-\{\ft\}$ is contained in $q^{d(\ft)}+1$ chambers \cite[\S 2.4]{tit}.
If $G$ is $k$-split (i.e. there is a maximal torus $T\subset G$ which is $k$-split) then $d(\ft)=1$ for all $\ft\in I$ \cite[\S 3.5.4]{tit}.

If $k$ has characteristic zero, then the condition that $\Gamma$ is torsion free
can be omitted from Theorem \ref{main}.
Recall that a non-archimedean local field of characteristic zero is a finite extension of $\QQ_p$,
for some prime $p$.  

\begin{corollary}\label{maincorollary} Let $k$ be a non-archimedean local field of characteristic zero. Let $G$ be the group of $k$-rational points of an absolutely almost simple, simply connected linear algebraic $k$-group.
If $\Gamma$ is a lattice in $G$, then the class $[\1]$ has torsion in $\Omega_\Gamma$.
\end{corollary}

\begin{proof}
 A lattice $\Gamma$ in $G$ is automatically cocompact \cite[Proposition IX, 3.7]{M}.  By Selberg's Lemma \cite[Theorem 2.7]{garland}, $\Gamma$ has a torsion free subgroup $\Gamma_0$ of finite index. Now Theorem \ref{main} implies that $[\1]$ has finite order in $\Omega_{\Gamma_0}$. 
The result follows from the observation that there is a natural surjection $\Omega_{\Gamma_0} \to \Omega_\Gamma$.
\end{proof}

\section{Proof of Theorem \ref{main}}\label{two}

Throughout this section, the assumptions of Theorem \ref{main} are in force.
Before proving Theorem \ref{main}, we require some preliminaries.
Recall that a gallery of type $i=(i_1,\ldots,i_k)$
is a sequence of chambers $(\sigma_0, \sigma_1,\ldots\sigma_k)$ such that each pair of successive
chambers $\sigma_{j-1}, \sigma_j$ meet in a common face of type $I-\{i_j\}$.  
Choose a special type $\fs\in I$, which will remain fixed throughout this section.
Fix once and for all the following data.
\begin{itemize}
\item{({\bf A1})} An apartment $A$ in $\Delta$.
\item{({\bf A2})} A sector $S$ in $A$ with base vertex $v$ of type $\fs$ and base chamber
            $C$.
\item{({\bf A3})} The unique vertex $v'\in S$ of type $\fs$, obtained by reflecting $v$ in a codimension-1 face of $C$.            
\item{({\bf A4})} The unique chamber $C'$ containing $v'$ which is the base chamber of a subsector of $S$.
\item{({\bf A5})} A minimal gallery of type $i=(i_1,\ldots,i_k)$ from $C$ to $C'$, where $i_1=\fs$.  This
        minimal gallery necessarily lies inside $S$.
\end{itemize}

\medskip

\begin{figure}[htbp]
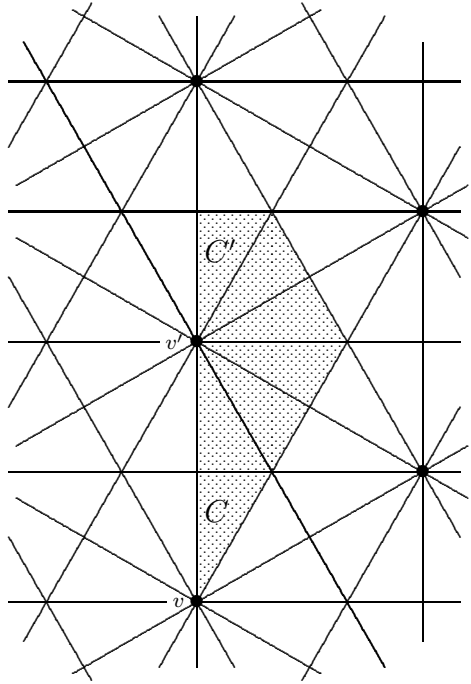
\label{apt1}
\centerline{
\beginpicture
\setcoordinatesystem units <1cm, 1.732cm>  
\setplotarea  x from -6 to 6,  y from -2.5 to 2.6
\putrule from -3.5 -2 to -1.4 -2  \putrule from -1 -2 to 2.5 -2  
\putrule from -3.5 -1 to 2.5 -1  
\putrule from -3.5 0 to -1.5 0     \putrule from -1 0 to 2.5 0
\putrule from -3.5 1  to 2.5  1
\putrule from -3.5 2  to 2.5  2
\putrule from  -1 -2.5  to -1 2.5
\putrule from  2 -2.3  to 2 2.3
\put {$C$} [l]     at   -0.9 -1.3
\put {$C^{\prime}$} [l]     at   -0.9 0.7
\put {$_v$}[r] at -1.15 -2 
\put {$_{v'}$}[r] at -1.15 0 
\put {$\bullet$} at -1 -2 
\put {$\bullet$} at -1  0  
\put {$\bullet$} at -1  2 
\put {$\bullet$} at  2  1  
\put {$\bullet$} at  2 -1
\setlinear
\plot -3.4   -0.8   -1 0      2.6    1.2   /
\plot -2.65  -2.55    -1 -2     2.6   -0.8  /
\plot  -3.4      1.2    -1  2      0.65    2.55   /
\plot   -3.4      0.8    -1 0      2.6   -1.2   /
\plot   -3.4     -1.2    -1 -2     0.65   -2.55   /
\plot -2.65      2.55   -1 2    2.6    0.8    /
\plot -3.5   1.5     -2.4  2.6  /
\plot -3.5  -0.5     -0.5  2.5  /
\plot -3.3  -2.3      1.5  2.5  /
\plot -1.5  -2.5      2.5  1.5  /
\plot  0.4  -2.6       2.5 -0.5  /
\plot -3.5  -1.5     -2.4 -2.6  /
\plot -3.5   0.5     -0.5 -2.5  /
\plot -3.3   2.3      1.5 -2.5  /
\plot -3.3   2.3      1.5 -2.5  /
\plot -1.5   2.5      2.5 -1.5  /
\plot  0.4   2.6      2.5  0.5  /
\setshadegrid span <1.5pt>
\vshade   -1 -2 1  <,z,,>  0 -1 1 <z,,,>  1 0 0 /
\endpicture
}
\hfil
\caption{Part of an apartment $A$ in a building of type $\widetilde G_2$.}
\end{figure}

These data are illustrated by Figure \ref{apt1}, which shows part of an apartment in a building of type $\widetilde G_2$ and a minimal gallery from $C$ to $C^{\prime}$.
 Special vertices are indicated by large points.

\refstepcounter{picture}

Now let $\fD=\Delta^n/\Gamma$, the set of $\Gamma$-orbits of chambers of 
$\Delta$.  Since $\Gamma$ acts freely and cocompactly on $\Delta$, $\fD$ is finite and elements of $\fD$ are in $1-1$ correspondence with the set of $n$-cells of the finite cell complex $\Delta/\Gamma$.

If $x$, $y\in\fD$, let $M_i(x,y)$ denote the number of $\Gamma$-orbits of galleries of type $i$
which have initial chamber in $x$ and final chamber in $y$.  
If $\sigma_0\in x$ is fixed then  $M_i(x,y)$ is equal to the number of galleries
$(\sigma_0, \sigma_1,\ldots\sigma_k)$ of type $i$ with final chamber $\sigma_k\in y$.
To see this, note that any gallery of type $i$ with initial chamber in $x$ and final chamber in $y$
lies in the $\Gamma$-orbit of such a gallery $(\sigma_0, \sigma_1,\ldots\sigma_k)$.
Moreover, two distinct galleries of this form lie in different $\Gamma$-orbits. 
For suppose that $(\sigma_0, \sigma_1,\ldots, \sigma_j, \tau_{j+1},\ldots \tau_k)$ is another such gallery,
with $\tau_{j+1}\ne\sigma_{j+1}$, the first chamber at which they differ. Then $\tau_{j+1}$ and $\sigma_{j+1}$ have a common face of codimension one, and so lie in different $\Gamma$-orbits, since the action of $\Gamma$ is free. (If $g\tau_{j+1}=\sigma_{j+1}$, then $g$ must fix every point in the common codimension one face
and so $g=1$.)   A similar argument shows that if $\sigma_k\in y$ is fixed then
then  $M_i(x,y)$ is equal to the number of galleries
$(\sigma_0, \sigma_1,\ldots\sigma_k)$ of type $i$ with initial chamber $\sigma_0\in x$.
Cocompactness of the $\Gamma$-action implies that $M_i(x,y)$ is finite.

If $\sigma$ is a chamber in $\Delta$, then the number $N_i$ of galleries 
$(\sigma_0, \sigma_1,\ldots,\sigma_k)$ of type $i$, with final chamber 
$\sigma_k=\sigma$, is independent of $\sigma$.  This follows, since $G$ acts transitively on the set $\Delta^n$ of chambers of $\Delta$.
Note that $N_i>1$, by thickness of the building $\Delta$.

Two different galleries of type $i$ which have final chamber $\sigma$ are necessarily in different $\Gamma$-orbits, by freeness of the action of $\Gamma$. It follows that if $y\in\fD$, then the number of $\Gamma$-orbits of galleries $(\sigma_1,\ldots,\sigma_k)$ of type $i$, with $\sigma_k\in y$, is equal to $N_i$.  In other words, for each $y\in\fD$, 
\begin{equation}\label{gallery cardinality}
\sum_{x\in\fD}M_i(x,y)=N_i\,.
\end{equation}

Recall that if $\tau$ is a simplex of $\Delta$ of codimension one and type $I-\{\ft\}$, then the number of chambers 
of $\Delta$ which contain $\tau$ is $q_{\ft}+1$ where $q_\ft\ge2$.
Thus the number of galleries $(\sigma_0, \sigma_1,\ldots,\sigma_k)$ of type $i$, with final chamber 
$\sigma_k=\sigma$ (fixed, but arbitrary), is equal to $q_{i_k}q_{i_{k-1}}\ldots q_{i_1}$,
where $i_1=\fs$.
On the other hand, this number is also equal to the number $q_{i_1}q_{i_2}\ldots q_{i_k}$ of
galleries $(\sigma_0, \sigma_1,\ldots,\sigma_k)$ of type $i$, with initial chamber 
$\sigma_0=\sigma$ (fixed, but arbitrary). It follows that, for each $x\in\fD$, 
\begin{equation}\label{gallery cardinality transposed}
\sum_{y\in\fD}M_i(x,y)=N_i\,.
\end{equation}

\begin{definition} Fix a type $\fs\in I$.
Let $\alpha_{\fs}$ denote the number of chambers of $\Delta$ which contain a fixed vertex $u$
of type $\fs$. Since $G$ acts transitively on the set of vertices of type $\fs$,  $\alpha_{\fs}$
does not depend on the choice of the vertex $u$.
\end{definition}

\begin{remark}\label{parahoricvolume}
The Iwahori subgroup $\mathfrak I$ is a chamber of $\Delta$.
Let the parahoric subgroup $P_\fs < G$ be the vertex of the type $\fs$ of $\mathfrak I$.
Then $P_\fs$ is a maximal compact subgroup of $G$ containing $\mathfrak I$ and $\alpha_\fs=|P_\fs : \mathfrak I|$.
(In \cite[Section 3]{garland}, $\alpha_\fs$ is denoted $\tau_{\{\fs\}}$.) 
\end{remark}

\begin{lemma}\label{Lone}
Let $\fs\in I$ be a special type.
Then 
\begin{equation}\label{estimate}
N_i< q_\fs\cdot\alpha_\fs.
\end{equation} 
\end{lemma}

\begin{figure}[htbp]
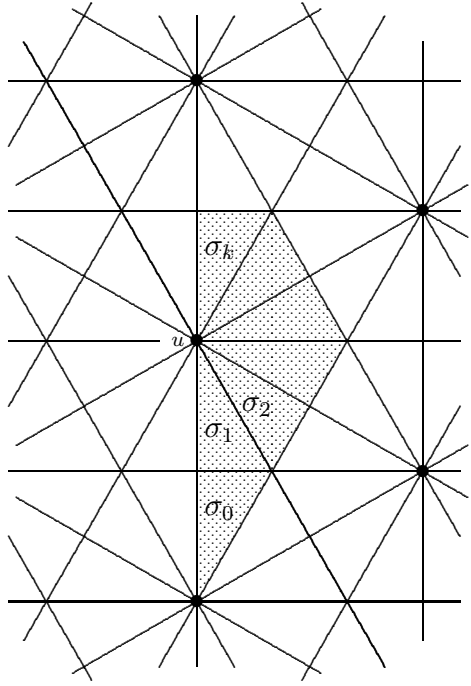
\label{gallery}
\centerline{
\beginpicture
\setcoordinatesystem units <1cm, 1.732cm>  
\setplotarea  x from -6 to 6,  y from -2.5 to 3
\putrule from -3.5 -2 to  2.5 -2
\putrule from -3.5 -1 to 2.5 -1  
\putrule from -3.5 0 to -1.5 0     \putrule from -1 0 to 2.5 0
\putrule from -3.5 1  to 2.5  1
\putrule from -3.5 2  to 2.5  2
\putrule from  -1 -2.5  to -1 2.5
\putrule from  2 -2.3  to 2 2.3
\put {$_{u}$}[r] at -1.15 0 
\put {$\bullet$} at -1 -2 
\put {$\bullet$} at -1  0  
\put {$\bullet$} at -1  2 
\put {$\bullet$} at  2  1  
\put {$\bullet$} at  2 -1
\put {$\sigma_0$} [l]     at   -0.9 -1.3
\put {$\sigma_1$} [l]     at   -0.9 -0.7
\put {$\sigma_2$} [r]     at   0 -0.5
\put {$\sigma_k$} [l]     at   -0.9 .7
\setlinear
\plot -3.4   -0.8   -1 0      2.6    1.2   /
\plot -2.65  -2.55    -1 -2     2.6   -0.8  /
\plot  -3.4      1.2    -1  2      0.65    2.55   /
\plot   -3.4      0.8    -1 0      2.6   -1.2   /
\plot   -3.4     -1.2    -1 -2     0.65   -2.55   /
\plot -2.65      2.55   -1 2    2.6    0.8    /
\plot -3.5   1.5     -2.4  2.6  /
\plot -3.5  -0.5     -0.5  2.5  /
\plot -3.3  -2.3      1.5  2.5  /
\plot -1.5  -2.5      2.5  1.5  /
\plot  0.4  -2.6       2.5 -0.5  /
\plot -3.5  -1.5     -2.4 -2.6  /
\plot -3.5   0.5     -0.5 -2.5  /
\plot -3.3   2.3      1.5 -2.5  /
\plot -3.3   2.3      1.5 -2.5  /
\plot -1.5   2.5      2.5 -1.5  /
\plot  0.4   2.6      2.5  0.5  /
\setshadegrid span <1.5pt>
\vshade   -1 -2 1  <,z,,>  0 -1 1 <z,,,>  1 0 0 /
\endpicture
}
\hfil
\caption{A minimal gallery $(\sigma_0, \sigma_1,\ldots,\sigma_k)$
in a $\widetilde G_2$ building.}
\end{figure}

\begin{proof}
Fix a chamber $\sigma_0$. We must estimate the number of galleries 
$(\sigma_0, \sigma_1,\ldots,\sigma_k)$ of type $i$ (with initial chamber $\sigma_0$).

There are $q_\fs$ possible choices for $\sigma_1$. Suppose that $\sigma_1$ has been chosen
and let $u$ be the vertex of $\sigma_1$ not belonging to $\sigma_0$.
By construction, $\sigma_k$ also contains $u$ (Figure \ref{gallery}) and so there are less than 
$\alpha_\fs$ possible choices for $\sigma_k$. (Note the $\sigma_k\ne\sigma_1$.) 
Once $\sigma_k$ has been chosen, there is a unique (minimal) gallery of type
$(i_2, \dots , i_k)$ with initial chamber $\sigma_1$ and final chamber $\sigma_k$.
In other words, the gallery $(\sigma_0, \sigma_1,\ldots,\sigma_k)$ is uniquely determined,
once $\sigma_1$ and $\sigma_k$ are chosen. There are therefore at most $q_\fs(\alpha_\fs-1)$
choices for this gallery.
\end{proof}

\begin{remark}
An easy calculation in $\widetilde A_2$ buildings shows that the estimate (\ref{estimate})
 cannot be improved to $N_i\le \alpha_\fs$.
\end{remark}

\begin{definition} Let $\Gamma$ be a torsion free cocompact lattice in $G$. 
If  $\fs\in I$, let $n_\fs(\Gamma)$ (or simply $n_\fs$, if $\Gamma$ is understood) denote the number of 
$\Gamma$-orbits of vertices of type $\fs$  in $\Delta$.
\end{definition}

\noindent Recall that $\covol(\Gamma)$ is equal to the number of $\Gamma$-orbits of chambers  in $\Delta$.

\begin{lemma}\label{Ltwo}
Fix a type $\fs\in I$.
Then $\covol(\Gamma)=n_\fs(\Gamma)\cdot \alpha_{\fs}$.  
\end{lemma}

\begin{proof}
Choose a set $\cS$ of representative vertices from the $\Gamma$-orbits of vertices of type $\fs$ in $\Delta$. Thus $|\cS|=n_\fs(\Gamma)$. 
For $v\in \cS$, let $R_v$ denote the set of chambers containing $v$. Each $R_v$ contains
$\alpha_{\fs}$ chambers. 
We claim that the number of chambers in $R=\displaystyle \bigcup_{v\in \cS} R_v$ equals $\covol(\Gamma)$.

Each chamber in $\Delta$ is clearly in the $\Gamma$-orbit of some chamber in $R$.  Moreover, any two distinct
chambers in $R$ lie in different $\Gamma$-orbits. For suppose that $\sigma_v\in R_v$, $\sigma_w\in R_w$ and
$g\sigma_v=\sigma_w$, where $g\in\Gamma$. Then $gv=w$, since the action of $\Gamma$ is type preserving
and every chamber contains exactly one vertex of type $\fs$. Therefore $v=w$, since distinct vertices in $\cS$ 
lie in different $\Gamma$-orbits.  Moreover $g=1$, since the action of $\Gamma$ is free.  Thus $\sigma_v=\sigma_w$.  This shows that there are $\covol(\Gamma)$ chambers in $R$.
\end{proof}

\refstepcounter{picture}
 
\medskip

Before proving Theorem \ref{main}, we provide more details of the structure of the boundary $\Omega$.
Let $\sigma$ be a chamber in $\Delta^n$ and let $s$ be a special vertex
of $\sigma$.  The codimension one faces of $\sigma$ having $s$ as a vertex
determine roots containing $\sigma$, and the intersection of these roots is a sector in $\Delta$ with base vertex $s$ and base chamber $\sigma$.
Two sectors are {\it parallel} if the Hausdorff distance
between them is finite.  This happens if and only if they contain a common subsector.
The boundary $\Omega$ of $\Delta$ is the set of parallel equivalence classes of sectors in $\Delta$ \cite[Chap. 9.3]{Ronan}.  If $\omega\in\Omega$ and if $s$ is a special vertex of $\Delta$ then there exists a unique sector $[s,\omega)$ in $\omega$ with base vertex $s$, \cite[Lemma 9.7]{Ronan}.

If $\sigma\in\Delta^n$, let $o(\sigma)$ denote the vertex of $\sigma$ of type $\fs$.  Recall that vertices of type $\fs$ are special.  Let $\Omega(\sigma)$ denote the set of boundary points $\omega$ whose representative sectors have base vertex $o(\sigma)$ and base chamber $\sigma$.
That is, 
\begin{equation*}
\Omega(\sigma)=\{\omega\in\Omega:\sigma\subset[o(\sigma),\omega)\}\,.
\end{equation*}
The sets $\Omega(\sigma)$, $\sigma\in\Delta^n$, form a basis for the topology of $\Omega$.
Moreover, each $\Omega(\sigma)$ is a clopen subset of $\Omega$.
Let $\gamma_i$ denote the set of ordered pairs $(\sigma,\sigma')\in\Delta^n\times\Delta^n$ such that there exists a gallery
of type $i$ from $\sigma$ to $\sigma'$.  Then for each $\sigma\in\Delta^n$, $\Omega(\sigma)$ can be expressed as a disjoint union
\begin{equation}\label{boundary sets}
\Omega(\sigma)=\bigsqcup_{(\sigma,\sigma')\in\gamma_i}\Omega(\sigma')\,.
\end{equation}
For if $\omega\in\Omega(\sigma)$, then the sector $[o(\sigma),\omega)$ is strongly isometric, in the sense of \cite[15.5]{garrett} to the sector
$S$ in the apartment $A$, as described at the beginning of this section.  Let $\sigma'$ be the image under this
strong isometry of the chamber $C'$ in $A$.  Then $(\sigma,\sigma')\in\gamma_i$ and  $\omega\in\Omega(\sigma')$.  Thus $\Omega(\sigma)$ is indeed a subset of the right hand side of 
(\ref{boundary sets}).  Conversely, each set $\Omega(\sigma')$ on the right hand side of (\ref{boundary sets}) is contained in $\Omega(\sigma)$.  For if $(\sigma,\sigma')\in\gamma_i$ and  $\omega\in\Omega(\sigma')$ then the strong isometry from $[o(\sigma'),\omega)$ onto $S^{\prime}$
extends to a strong isometry from $[o(\sigma),\omega)$ onto $S$ \cite[\S 15.5 Lemma]{garrett}. Thus $\omega\in \Omega(\sigma)$.

To check that the union on the right of (\ref{boundary sets}) is disjoint, suppose that $\omega\in\Omega(\sigma'_1)\cap\Omega(\sigma'_2)$, where $(\sigma,\sigma'_1), (\sigma,\sigma'_2)\in\gamma_i$.  Then the strong isometry from $[o(\sigma'_1),\omega)$ onto $[o(\sigma'_2),\omega)$ extends to a strong isometry from $[o(\sigma),\omega)$ onto itself, which is necessarily the identity map. In particular, $\sigma'_1=\sigma'_2$.

If $\sigma\in\Delta^n$, let $\chi_\sigma\in C(\Omega,\ZZ)$ denote the characteristic function of $\Omega(\sigma)$.  That is
$$\chi_\sigma(\omega)=\begin{cases} 1\qquad\text{if $\omega\in\Omega(\sigma)$},\\
                        0\qquad\text{otherwise.}
               \end{cases}
$$
Since $\chi_\sigma-\chi_{g\sigma}=\chi_\sigma-g.\chi_\sigma$ for each $g\in\Gamma$, the class
$\left[\chi_\sigma\right]$ of $\chi_\sigma$ in $\Omega_\Gamma$ depends only on the $\Gamma$-orbit
of $\sigma$ in $\Delta^n$.  If $x=\Gamma\sigma\in\fD$, it therefore makes sense to define
\begin{equation}\label{class of element}
[x]=\left[\chi_\sigma\right]\in \Omega_\Gamma\,.
\end{equation}
Now it follows from (\ref{boundary sets}) that, for each $\sigma\in\Delta^n$, 
\begin{equation}\label{careful}
\chi_\sigma=\sum_{(\sigma,\sigma')\in\gamma_i}\chi_{\sigma'}
   =\sum_{y\in\fD}\sum_{  \begin{smallmatrix}\sigma'\in y\\(\sigma, \sigma')\in\gamma_i\end{smallmatrix}}\chi_{\sigma'}\,.
\end{equation}
Passing to equivalence classes in $\Omega_\Gamma$ gives, for each $x\in\fD$,
\begin{equation}\label{class2 of element}
[x]=\sum_{y\in\fD}M_i(x,y)[y]\,.
\end{equation}

\bigskip
We can now proceed with the proof of the Theorem \ref{main}.
If $s$ is a vertex of type $\fs$ of $\Delta$, then each element $\omega\in\Omega$ lies in $\Omega(\sigma)$ where $\sigma$ is the base chamber of the sector $[s,\omega)$.  Moreover $\omega$ lies in precisely one such set $\Omega(\sigma)$, with $\sigma\in\Delta^n$,
$o(\sigma)=s$.  Therefore
\begin{equation}\label{identity sum}
\1=\sum_{  \begin{smallmatrix}\sigma\in\Delta^n\\o(\sigma)=s\end{smallmatrix}}\chi_\sigma\,.
\end{equation}
Since the action of $\Gamma$ on $\Delta$ is free and type preserving, no two chambers $\sigma\in\Delta^n$ with $o(\sigma)=s$ lie in the same $\Gamma$-orbit.  To simplify notation, let $n_\fs=n_\fs(\Gamma)$, the number of $\Gamma$-orbits of vertices of type $\fs$ in $\Delta$.  If we choose a representative set $\cS$ of vertices of type $\fs$ in $\Delta$ then the chambers containing these vertices form a representative set of chambers, by the proof of Lemma \ref{Ltwo}.  It follows that in $\Omega_\Gamma$, 
\begin{alignat*}{2}
n_\fs\cdot[\1]
&=\sum_{s\in \cS}\sum_{\begin{smallmatrix}\sigma\in\Delta^n\\
                    o(\sigma)=s
                    \end{smallmatrix}}[\chi_\sigma]
            & &{ }\qquad\text{(by (\ref{identity sum}))}\\
&=\sum_{x\in\fD}[x] \,. & &
\end{alignat*}
Therefore
\begin{alignat*}{2}
n_\fs\cdot[\1]
&=\sum_{x\in\fD}\sum_{y\in\fD}M_i(x,y)[y]& &{ }\qquad
    \text{(by (\ref{class2 of element}))}\\
&=\sum_{y\in\fD}\left(\sum_{x\in\fD}M_i(x,y)\right)[y] & &\\
&=\sum_{y\in\fD}N_i\cdot [y]& &{ }\qquad\text{(by (\ref{gallery cardinality}))}\\
&= N_i n_\fs \cdot[\1]\,. & &
\end{alignat*}
It follows that 
\begin{equation}\label{ordereq}
n_\fs(N_i-1)\cdot[\1]=0,
\end{equation}
which proves the first assertion of Theorem \ref{main}.

Using Lemmas \ref{Lone}, \ref{Ltwo}, we can estimate the order of the element $[\1]$.
\begin{equation}
n_\fs(N_i-1) < n_\fs\cdot(q_\fs\alpha_\fs-1) =q_\fs\cdot \covol(\Gamma)-n_\fs.
\end{equation}
This proves the second assertion of Theorem \ref{main}.
The next Lemma proves the final assertion of Theorem \ref{main} by showing that the estimate of the order of $[\1]$ can be improved if certain exceptional cases are excluded.

\begin{lemma}\label{improvedestimate}
Suppose that the Weyl group is not one of the exceptional types $\widetilde E_8, \widetilde F_4, \widetilde G_2$. Then
\[
\ord([\1]) < \covol(\Gamma) \, .
\]
\end{lemma}

\begin{proof}
An examination of the possible Coxeter diagrams \cite[Chap VI, No 4.4, Th\'eor\`eme 4]{bou}
shows that if the diagram is not one of the types $\widetilde E_8, \widetilde F_4, \widetilde G_2$,
then it contains at least two special types. 
Therefore every chamber of $\Delta$ contains at least two special vertices. Choose two such vertices
and suppose that they have types $\fs$ and $\ft$, say.  
In that case the condition {\bf (A3)} on  the apartment $A$ in $\Delta$ can be changed to read:
\begin{itemize}
\item{({\bf A3$^\prime$})} The unique special vertex $v'\in S$ of type $\ft$ which lies in $C$.            
\end{itemize}
Assume that the remaining conditions {\bf (A1)}, {\bf (A2)}, {\bf (A4)}, {\bf (A5)}
are unchanged.
Figure \ref{aptB2} illustrates the setup in the $\widetilde B_2$ case.

\refstepcounter{picture}
\begin{figure}[htbp]
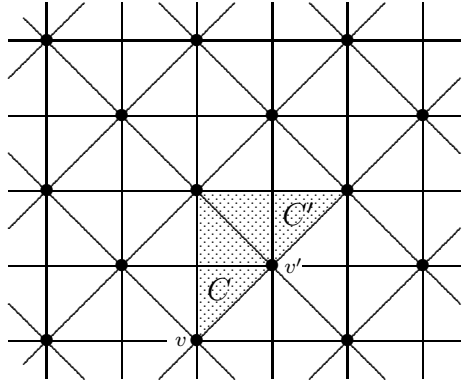
\label{aptB2}
\centerline{
\beginpicture
\setcoordinatesystem units <1cm, 1cm>  
\setplotarea  x from -6 to 6,  y from -2.5 to 3
\putrule from -3.5 -2 to -1.4 -2  \putrule from -1 -2 to 2.5 -2   
\putrule from -3.5 -1 to  0 -1  \putrule from 0.4  -1 to 2.5 -1  
\putrule from -3.5 0 to  2.5 0 
\putrule from -3.5 1  to 2.5  1
\putrule from -3.5 2  to 2.5  2
\putrule from  -3 -2.5  to -3 2.5
\putrule from  -2 -2.5  to -2 2.5
\putrule from  -1 -2.5  to -1 2.5
\putrule from  0 -2.5  to -0 2.5
\putrule from  1 -2.5  to 1 2.5
\putrule from  2 -2.5  to 2 2.5
\put {$_v$}     at   -1.2 -2
\put {$_{v'}$}   at   0.3  -1
\put {$C$}     at   -0.7 -1.3
\put {$C^{\prime}$}     at   0.35  -0.3
\multiput {$\bullet$} at -3 -2 *2  2 0 /
\multiput {$\bullet$} at -2 -1 *2  2 0 /
\multiput {$\bullet$} at -3 0 *2  2 0 /
\multiput {$\bullet$} at -2 1 *2  2 0 /
\multiput {$\bullet$} at -3 2 *2  2 0 /
\setlinear
\plot -3.5 1.5     -2.5 2.5 /
\plot -3.5 -0.5     -0.5 2.5 /
\plot -3.3 -2.3     1.5 2.5 /
\plot -1.5 -2.5     2.5 1.5 /
\plot 0.5 -2.5     2.5 -0.5 /
\plot -3.5 -1.5     -2.5 -2.5 /
\plot -3.5 .5     -0.5 -2.5 /
\plot -3.3 2.3     1.5 -2.5 /
\plot -1.5 2.5     2.5 -1.5 /
\plot 0.5 2.5     2.5 0.5 /
\setshadegrid span <1.5pt>
\vshade   -1 -2 0  <z,,,>  1 0 0 /
\endpicture
}
\hfil
\caption{Part of an apartment $A$ in a building of type $\widetilde B_2$,
and a minimal gallery from $C$ to $C^{\prime}$.}
\end{figure}

The proof proceeds exactly as before, except that all the chambers in a gallery
$(\sigma_0, \sigma_1,\ldots,\sigma_k)$ of type $i$ now contain a common vertex 
$u$ of type $\ft$. Therefore equation(\ref{estimate}) becomes 
\begin{equation*}\label{newestimate}
N_i< \alpha_\ft\, .
\end{equation*}
Observe that one must be careful with the notation. For example in equation (\ref{careful}),
the function $\chi_\sigma$ on the left is now defined in terms of sectors based at the vertex 
of type $\fs$ of $\sigma$, whereas the functions $\chi_{\sigma'}$ on the right will now be defined in terms of sectors based at the vertex 
of type $\ft$ of $\sigma'$.
Equation (\ref{ordereq}) becomes
\begin{equation}\label{ordereqnew}
(n_\ft N_i- n_\fs)\cdot[\1]=0\, .
\end{equation} 
The order of the element $[\1]$ is bounded by
\begin{equation}\label{ten}
n_\ft\cdot\alpha_\ft-n_\fs  = \covol(\Gamma)-n_\fs.
\end{equation}
\end{proof}

Finally, we verify that $\Omega_\Gamma$ is a finitely generated group.
Sets of the form $\Omega(\sigma)$, $\sigma\in\Delta^n$, form a basis of clopen sets for the topology of $\Omega$.  It follows that the abelian group $C(\Omega,\ZZ)$ is generated by the set of characteristic functions $\{\chi_\sigma:\sigma\in\Delta^n\}$.  We show that $\Omega_\Gamma$ is generated by $\{[x]:x\in\fD\}$.

\begin{lemma}
Every clopen set $V$ in $\Omega$ may be expressed as a finite disjoint union of sets of the form $\Omega(\sigma)$, $\sigma\in\Delta^n$.
\end{lemma}

\begin{proof}
Fix a special vertex $s$ of type $\fs$ in $\Delta$. For each $\omega\in\Omega$, sets of the form $\Omega(\sigma)$ with $\sigma\in\Delta^n$ and $\sigma\subset [s,\omega)$ form a basic family of open neighbourhoods of $\omega$.  Therefore, for each $\omega\in V$, there exists a chamber $\sigma_\omega\in\Delta^n$ with $\sigma_\omega\subset[s,\omega)$ and $\omega\in\Omega(\sigma_\omega)\subseteq V$.  
The clopen set $V$, being compact, is a finite union of such sets:
$$V=\Omega(\sigma_{\omega_1})\cup\cdots\cup\Omega(\sigma_{\omega_k})\,.$$

Fix a sector $Q$ in $\Delta$, with base vertex $s$.  For each $j$, $1\leq j\leq k$, let $C_j$ be the chamber in $Q$ which is the image of $\sigma_{\omega_j}$, under the unique strong isometry from $[s,\omega_j)$ onto $Q$.  Let $Q_j$ be the subsector of $Q$ with base chamber $C_j$ ($1\leq j\leq k)$, and choose a chamber $C$ in $\bigcap^k_{j=1}Q_j$. Informally,
$C$ is chosen to be sufficiently far away from the base vertex $s$.

For $1\leq j\leq k$, let $\tau_j$ be the chamber in $[s,\omega_j)$ which is the image of $C$ under the strong isometry from $Q$ onto $[s,\omega_j)$.  For each $\omega\in\Omega(\sigma_{\omega_j})$ there is a retraction from $[s,\omega)$ onto $[s,\omega_j)$ \cite[4.2]{garrett}.
Let $\tau_j(\omega)$ be the inverse image of the chamber $\tau_j$ under this retraction.  By local finiteness of $\Delta$, there are only finitely many such chambers $\tau_j(\omega)$, $\omega\in\Omega(\sigma_{\omega_j})$.  Call them $\tau_{j,l}$, $1\leq l\leq n_j$.  Thus $\Omega(\sigma_{\omega_j})$ may be expressed as a finite disjoint union:
$$\Omega(\sigma_j)=\bigsqcup_l\Omega(\tau_{j,l})\,.$$
Moreover, if $\omega\in\Omega(\tau_{j,l})$ then the strong isometry from $[s,\omega)$ onto $Q$ maps $\tau_{j,l}$ to the chamber $C$.
Finally, $V$ may be expressed as a disjoint union: 
$$V=\bigsqcup_{j,l}\Omega(\tau_{j,l})\,.$$
To check that this union is indeed  disjoint, suppose that $\omega\in\Omega(\tau_{j,l})\cap\Omega(\tau_{r,s})$.  Then, under the strong isometry from $Q$ onto $[s,\omega)$, the image of the chamber $C$ is equal to both $\tau_{j,l}$ and $\tau_{r,s}$.  In particular, $\tau_{j,l}=\tau_{r,s}$.
\end{proof}

\begin{proposition}\label{characterization} 
Let $(G,{\mathfrak I},N,S)$ be an affine topological Tits system, and let $\Gamma$ be a subgroup
of $G$. Then
\begin{itemize}
\item[{\it (a)}] The abelian group $C(\Omega,\ZZ)$ is generated by the set of characteristic functions $\{\chi_\sigma:\sigma\in\Delta^n\}$.
\item[{\it (b)}] $\Omega_\Gamma$ is generated by $\{[x]:x\in\fD\}$. 
\end{itemize}
\end{proposition}

\begin{proof}
(a) Any function $f\in C(\Omega,\ZZ)$ is bounded, by compactness of $\Omega$, and so takes finitely many values $n_i\in\ZZ$.
Now $V_i=\{\omega\in\Omega:f(\omega)=n_i\}$ is a clopen set in $\Omega$.  It follows from the preceding Lemma that $f$ may be expressed as a finite sum $f=\sum_j m_j\chi_{\sigma_j}$, with
$\sigma_j\in\Delta^n$.

(b)  This is an immediate consequence of (a).
\end{proof}

\vskip 1cm

\section{Further calculations in the rank 2 case}\label{rank2}

This section is devoted to showing that the estimate for the order of $[\1]$ given by Theorem \ref{main} can be improved if the building $\Delta$ is 2-dimensional.
 The group $G$ has type $\widetilde A_2$, $\widetilde B_2$ or  $\widetilde G_2$.
Denote the type set by $I=\{\fs, \fa, \fb\}$, where $\fs$ is a special type of the corresponding Coxeter diagram, as indicated below. 
Note that in the $\widetilde B_2$ case, the vertex $\fb$ is also special. In the $\widetilde A_2$ case, all vertices are special and $q_\ft=q$ for all $\ft\in I$.

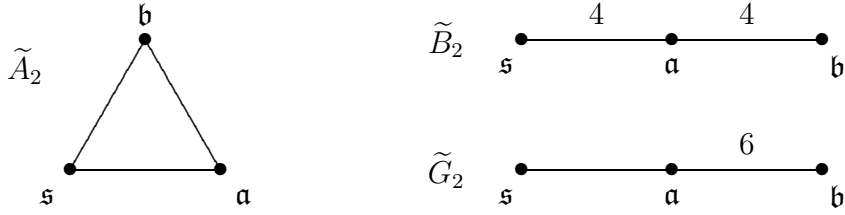
\centerline{
\beginpicture
\setcoordinatesystem units <1.0cm, 1.732cm>
\setplotarea x from -5 to 5, y from 0.3 to 2.5         
\putrule from -5 1 to -3 1  
\setlinear \plot -5 1   -4 2   -3 1 /
\put {$\bullet$} at -5 1 \put {$\bullet$} at -3 1 \put {$\bullet$} at -4 2
\put {$\fs$}  at -5.3 0.8  \put {$\fa$}  at -2.7 0.8  \put {$\fb$}  at -4 2.2
\put {$\widetilde A_2$}  at -5.6 1.8
\putrule from 1 2  to 5 2   
\put {$\bullet$} at 1 2 \put {$\bullet$} at 3 2 \put {$\bullet$} at 5 2
\put {$\fs$}  at 0.8 1.8  \put {$\fa$}  at 3 1.8  \put {$\fb$}  at 5.2 1.8
\put {$4$}  at 2 2.2  \put {$4$}  at  4 2.2
\put {$\widetilde B_2$}  at 0 2
\putrule from 1 1 to 5 1    
\put {$\bullet$} at 1 1 \put {$\bullet$} at 3 1 \put {$\bullet$} at 5 1
\put {$\fs$}  at 0.8 0.8  \put {$\fa$}  at 3 0.8  \put {$\fb$}  at 5.2 0.8
\put {$6$}  at 4 1.2
\put {$\widetilde G_2$}  at 0 1
\endpicture
}

\begin{proposition}\label{rank2prop}
Under the preceding assumptions, let $\Gamma$ be a torsion free lattice in $G$. Then
\begin{equation} \label{calc}
(q_\fa^2-1) n_\fs \cdot[\1] =  0
\end{equation}
in $\Omega_\Gamma$.
\end{proposition}

\begin{proof}
We prove the $\widetilde G_2$ case.
For the minimal gallery of type $i$ between $C$ and $C^\prime$ described in Figure \ref{apt1}, we obtain $N_i = q_\fs q_\fa^3q_\fb^2$, so that
\begin{equation} \label{calc1}
q_\fs q_\fa^3q_\fb^2 n_\fs \cdot[\1] =  n_\fs\cdot [\1].
\end{equation}
On the other hand, for a minimal gallery of type $j$ between $C$ and $C^{\prime\prime}$ described in Figure \ref{aptcalc1} below, we obtain $N_j = q_\fs^2 q_\fa^4q_\fb^4$, so that
\begin{equation} \label{calc2}
q_\fs^2 q_\fa^4q_\fb^4\ n_\fs \cdot[\1] =  n_\fs \cdot[\1].
\end{equation}
Equations (\ref{calc1}), (\ref{calc2}) imply that
$$q_\fa^2 n_\fs \cdot[\1] =  n_\fs \cdot[\1],$$
thereby proving (\ref{calc}).

The $\widetilde B_2$ and $\widetilde A_2$ cases follow by similar calculations, using the configurations in Figure \ref{aptcalc2} below.
\end{proof}

\refstepcounter{picture}
\begin{figure}[htbp]\label{aptcalc1}
\centerline{
\beginpicture
\setcoordinatesystem units <1cm, 1.732cm>  
\setplotarea  x from -6 to 6,  y from -2 to 2.5
\putrule from -3.5 -2 to -1.4 -2  \putrule from -1 -2 to 2.5 -2  
\putrule from -3.5 -1 to 2.5 -1  
\putrule from -3.5 0 to -1.5 0     \putrule from -1 0 to 2.5 0
\putrule from -3.5 1  to 1.6  1  \putrule from 2 1  to 2.5  1
\putrule from -3.5 2  to 3  2   
\putrule from  -1 -2.5  to -1 2.5
\putrule from  2 -2.3  to 2 2.3
\put {$C$} [l]     at   -0.9 -1.3
\put {$C^{\prime}$} [l]     at   -0.9 0.7
\put {$C^{\prime\prime}$} [l]     at   2.1 1.7
\put {$_\fs$}[r] at -1.15 -2 
\put {$_{\fs}$}[r] at -1.15 0 
\put {$_{\fs}$}[r] at 1.85 1 
\put {$_{\fa}$}  at -1.1 -0.9 
\put {$_{\fb}$} at 0.3 -0.9 
\put {$_{\fa}$} at 0.75   -0.45 
\put {$_{\fb}$} at 1.3   0.1
\put {$_{\fa}$} at  0.75  0.45 
\put {$_{\fb}$} at 0.3 1.1 
\put {$_{\fa}$} at  0.75  1.55
\put {$_{\fb}$} at 1.3   1.9
\put {$\bullet$} at -1 -2 
\put {$\bullet$} at -1  0  
\put {$\bullet$} at -1  2 
\put {$\bullet$} at  2  1  
\put {$\bullet$} at  2 -1
\setlinear
\plot -3.4   -0.8   -1 0      2.6    1.2   /
\plot -2.65  -2.55    -1 -2     2.6   -0.8  /
\plot  -3.4      1.2    -1  2      0.65    2.55   /
\plot   -3.4      0.8    -1 0      2.6   -1.2   /
\plot   -3.4     -1.2    -1 -2     0.65   -2.55   /
\plot -2.65      2.55   -1 2    2.6    0.8    /
\plot -3.5   1.5     -2.4  2.6  /
\plot -3.5  -0.5     -0.5  2.5  /
\plot -3.3  -2.3      1.5  2.5  /
\plot -1.5  -2.5      3  2 / 
\plot  0.4  -2.6       2.5 -0.5  /
\plot -3.5  -1.5     -2.4 -2.6  /
\plot -3.5   0.5     -0.5 -2.5  /
\plot -3.3   2.3      1.5 -2.5  /
\plot -3.3   2.3      1.5 -2.5  /
\plot -1.5   2.5      2.5 -1.5  /
\plot  0.4   2.6      2.5  0.5  /
\setshadegrid span <1.5pt>
\vshade   -1 -2 -1  <,z,,>  0 -1 -1 /
\vshade   -1 0 1  <,z,,>  0 1 1  /
\vshade   2   1 2  <,z,,>  3 2 2 /
\endpicture
}
\hfil
\caption{The $\widetilde G_2$ case.}
\end{figure}

\refstepcounter{picture}
\begin{figure}[htbp]\label{aptcalc2}
\centerline{
\beginpicture
\setcoordinatesystem units <1cm, 1cm>  
\setplotarea  x from -4 to 10,  y from -4 to 4
\put{  
\beginpicture
\setcoordinatesystem units  <0.5cm, 0.866cm>        
\setplotarea x from -2.5 to 2.5, y from -1 to 2.5   
\put {The $\widetilde A_2$ case.}  at   1 -3 
\put {$_\fs$}     at   -2.4 -2
\put {$_\fs$}     at   -2.4  0
\put {$_\fs$}     at   1.4  1
\put {$C$}     at   -1.9 -1.3
\put {$C^{\prime}$}     at  -1.8  0.7
\put {$C^{\prime\prime}$}     at   1.2  1.7
\put {$\bullet$} at  -2  -2 
\multiput {$\bullet$} at -3 -1 *1  2  0 /
\multiput {$\bullet$} at -4 0 *2  2 0 /
\multiput {$\bullet$} at -3 1 *2  2 0 /
\multiput {$\bullet$} at -2 2 *2  2 0 /
\putrule from -2.5   2     to  2.5  2
\putrule from  -3.5 1  to  1  1  
\putrule from -4.5 0  to -2.6 0  \putrule from -2 0  to 0.5 0    
\putrule from -3.5 -1  to  -0.5 -1
\setlinear
\plot  -4.3 -0.3  -1.7 2.3 /
\plot -3.3 -1.3  0.3 2.3 /
\plot -2.2 -2.2  2.3 2.3 / 
\plot  -4.3 0.3  -1.8 -2.2 /
\plot  -3.3 1.3   -0.7 -1.3 / 
\plot  -2.3 2.3   0.3 -0.3 /
\plot  -0.3 2.3   1.3 0.7 /
\plot  1.7 2.3   2.3 1.7 /
\setshadegrid span <1.5pt>
\vshade   -3 -1 -1  <,z,,>  -2 -2 -1  <z,,,>  -1   -1 -1 /
\vshade    -3  1 1  <,z,,> -2  0 1  <z,,,>  -1   1  1 /
\vshade   0  2  2  <,z,,>  1  1 2  <z,,,>  2   2 2 /
\endpicture
}
at  6  -0.5   
\put {The $\widetilde B_2$ case.}     at   -0.5 -3
\putrule from -3.5 -2 to -1.4 -2  \putrule from -1 -2 to 2.5 -2   
\putrule from -3.5 -1 to  0 -1  \putrule from 0.3  -1 to 2.5 -1  
\putrule from -3.5 0 to  -1.4  0 \putrule from -1 0 to  0.6  0 \putrule from 1  0  to  2.5 0 
\putrule from -3.5 1  to 2.5  1
\putrule from -3.5 2  to 2.5  2
\putrule from  -3 -2.5  to -3 2.5
\putrule from  -2 -2.5  to -2 2.5
\putrule from  -1 -2.5  to -1 2.5
\putrule from  0 -2.5  to -0 2.5
\putrule from  1 -2.5  to 1 2.5
\putrule from  2 -2.5  to 2 2.5
\put {$_\fs$}     at   -1.2 -2
\put {$_\fs$}     at   -1.2  0
\put {$_\fs$}     at   0.8 0
\put {$_\fb$}   at   0.2  -1
\put {$_\fa$}   at   -1.15  -0.85
\put {$C$}     at   -0.7 -1.3
\put {$C^{\prime}$}     at  -0.7 0.7
\put {$C^{\prime\prime}$}     at   1.35  0.7
\multiput {$\bullet$} at -3 -2 *2  2 0 /
\multiput {$\bullet$} at -2 -1 *2  2 0 /
\multiput {$\bullet$} at -3 0 *2  2 0 /
\multiput {$\bullet$} at -2 1 *2  2 0 /
\multiput {$\bullet$} at -3 2 *2  2 0 /
\setlinear
\plot -3.5 1.5     -2.5 2.5 /
\plot -3.5 -0.5     -0.5 2.5 /
\plot -3.3 -2.3     1.5 2.5 /
\plot -1.5 -2.5     2.5 1.5 /
\plot 0.5 -2.5     2.5 -0.5 /
\plot -3.5 -1.5     -2.5 -2.5 /
\plot -3.5 .5     -0.5 -2.5 /
\plot -3.3 2.3     1.5 -2.5 /
\plot -1.5 2.5     2.5 -1.5 /
\plot 0.5 2.5     2.5 0.5 /
\setshadegrid span <1.5pt>
\vshade   -1 -2 -1  <z,,,>  0   -1 -1 /
\vshade   -1  0 1  <z,,,>  0   1  1 /
\vshade   1  0  1  <z,,,>  2   1  1 /
\endpicture
}
\hfil
\caption{}
\end{figure}

Let $k$ be a non-archimedean local field with residue field $\overline k$ of order $q$. 
Let $L$ be a simple, simply connected linear algebraic $k$-group and assume that $L$ is 
$k$-split and has $k$-rank 2. Let $G$ be the group of $k$-rational points of $L$ and let $\Gamma$ be a torsion free lattice in $G$.
Then $q_\ft=q$ for all $\ft\in I$ \cite[\S 3.5.4]{tit}, and equation (\ref{calc}) becomes
\begin{equation} \label{splitcalc}
(q^2-1) n_\fs \cdot[\1] =  0.
\end{equation}
A parahoric subgroup $P_\fs$ corresponding to a {\it hyperspecial} vertex of type $\fs$ has maximal volume among compact subgroups of $G$ \cite[3.8.2]{tit}. This volume is $[P_\fs : \mathfrak I]$, by Remark \ref{parahoricvolume}. In particular, all such subgroups have the same volume. 
It follows that $n_\fs=\covol(\Gamma)/[P_\fs : \mathfrak I]$ has the same value for all hyperspecial types $\fs$.

Suppose, for example, that $G$ is the symplectic group $Sp_2(k)$, which has type $\widetilde B_2$ (or, equivalently, $\widetilde C_2$). 
Examination of the tables at the end of \cite{tit} shows that the diagram of $G$ has two hyperspecial types $\fs$, $\ft$. 
Thus $n_\fs=n_\ft$, and it follows from (\ref{ordereqnew}) and Figure \ref{aptB2}, that 
\begin{equation*}
(q^3-1) n_\fs \cdot[\1] =  0.
\end{equation*}
Combining this with (\ref{splitcalc}) gives the following improvement to (\ref{splitcalc}),
for the case $G=Sp_2(k)$\,:
\begin{equation}\label{newsplitcalc}
(q-1) n_\fs \cdot[\1] =  0.
\end{equation}

\begin{remark}
An interesting problem is to find the exact value of the order of $[\1]$.
This is known in the case where the group $G$ has k-rank 1, and $\Delta$ is a tree. In that case a torsion free lattice $\Gamma$ in $G$ is a free group of finite rank $r$, and it follows from \cite {R, R2} that $[\1]$ has order $r-1=-\chi(\G)$, where $\chi(\Gamma)$ denotes
the Euler-Poincar\'e characteristic of $\Gamma$. If $G=\SL_2(k)$, then $-\chi(\Gamma)=(q-1)n_\fs(\Gamma)$.

In the rank 2 case, the order of $[\1]$ is in general smaller than $\chi(\Gamma)$.
For by \cite[p. 150, Th\'eor\`eme 7]{ser0}, $\chi(\Gamma)=(q-1)(q^m-1)n_\fs(\Gamma)$, where $m=2,3,5$ according as $G$ has type $\widetilde A_2$, $\widetilde B_2$, $\widetilde G_2$. Note however that by (\ref{splitcalc}), (\ref{newsplitcalc}), we do have $\chi(\Gamma)\cdot[\1]=0$
if $G=\SL_3(k)$ or $G=Sp_2(k)$.
\end{remark}

\vskip 1cm

\section{K-Theory of the Boundary Algebra $\cA_\Gamma$}\label{boundaryalgebra}

We retain the general assumptions of Theorem \ref{main}. Thus $G$ is a locally compact group acting strongly transitively by type preserving automorphisms on the affine building $\Delta$, and $\Gamma$ is a torsion free discrete subgroup of $G$.

As in \cite{RS}, \cite{R}, the group $\Gamma$ acts on the commutative $C^*$-algebra $C(\Omega)$, and one can form the full crossed product $C^*$-algebra $\cA_\Gamma=C(\Omega)\rtimes\Gamma$, \cite[Section 1]{R}.  The inclusion map $C(\Omega)\to \cA_\Gamma$ induces a natural homomorphism from $C(\Omega,\ZZ)=K_0(C(\Omega))$ to $K_0(\cA_\Gamma)$, which maps $\chi_\sigma$ to the class of the corresponding idempotent in $\cA_\Gamma$.  The covariance relations in $\cA_\Gamma$ imply that for each $g\in\Gamma$ and $\sigma\in\Delta^n$, the functions $\chi_\sigma$ and $g\cdot\chi_\sigma=\chi_{g\sigma}$ map to the same element of $K_0(\cA_\Gamma)$.  Thus there is an induced homomorphism $\varphi:\Omega_\Gamma\to K_0(\cA_\Gamma)$.  Moreover $\varphi([\1])=[\1]_{K_0}$, the class of \1 in the $K_0$-group of $\cA_\Gamma$. We have the following immediate consequence of Theorem \ref{main}.

\medskip

\begin{corollary}\label{abt}
If\, $\Gamma$ is a torsion free lattice in $G$ then $[\1]_{K_0}$ has finite order in $K_0(\cA_\Gamma)$.
\end{corollary}

\medskip

\begin{remark}
Clearly the bounds for the order of $[\1]$ obtained in the preceding sections apply also to $[\1]_{K_0}$.
If $G$ has type $\widetilde A_n$, Corollary \ref{abt} was proved in \cite{RS}, \cite{R}. In that case $q_\ft=q$ for all $\ft\in I$. For $n=1$, it follows from \cite {R, R2} that the order of $[\1]_{K_0}$ is actually
\begin{equation}\label{conjecture} 
\ord([\1]_{K_0})= (q-1)\cdot n_\fs.
\end{equation}
The computational evidence at the end of Section \ref{a2tilde} below indicates that (\ref{conjecture}) also holds for $n=2$. 
\end{remark}

Return now to the general assumptions of Theorem \ref{main}. It is important that $\Gamma$ is amenable at infinity in the sense of \cite[Section 5.2]{AR}. 
Since the action of $G$ on $\Delta$ is strongly transitive, its action on the boundary $\Omega$ is transitive. Therefore $\Omega$ may be identified, as a topological $\Gamma$-space, with $G/B$, where the Borel subgroup $B$ is the stabilizer of some point $\omega\in\Omega$.  The next result shows that the group $B$ is amenable and so the action of $\Gamma$ on $\Omega$ is amenable \cite[Section 2.2]{AR}.  Moreover the crossed product algebra $\cA_\Gamma$ is unique : the full and reduced crossed products coincide.

\begin{proposition}\label{amenable}
Let $\omega\in\Omega$ and let $B=\{g\in G:g\omega=\omega\}$.  Then $B$ is amenable and so $(\Gamma, \Omega)$ is amenable as a topological $\Gamma$-space, if $\Gamma$
is a closed subgroup of $G$.
\end{proposition}

\begin{proof}
Let $s\in\Delta^0$ be a special vertex and let $A$ be an apartment in $\Delta$ containing the sector $[s,\omega)$.  Let $N_{{\rm trans}}$ denote the subgroup of $G$ consisting of elements which stabilize $A$ and act by translation on $A$.

If $g\in B$, then the sectors $[gs,\omega)$ and $g[s,\omega)$ both have base vertex $gs$ and both represent the same boundary point $\omega$.  Therefore $g[s,\omega)=[gs,\omega)$.
Now the sectors $[gs,\omega)$, $[s,\omega)$ and $[g^{-1}s,\omega)$ are all equivalent, and so contain a common subsector $S$.  The sectors $S$ and $gS$, being subsectors of $[s,\omega)$, are parallel sectors in the apartment $A$.
Let $\sigma$ be the base chamber of $S$. Since $G$ acts strongly transitively on $\Delta$, there exists
an element $g'\in G$ such that $g' A= A$ and $g'\sigma
= g\sigma$. In particular $g'\omega = \omega$.

Since the action of $G$ is type preserving, it follows from \cite[Theorem 17.3]{garrett} that $g'\in N_{{\rm trans}}$.
Moreover $gv=g'v$, for all $v\in S$.
Let $\lambda_\omega(g)=g'|_A$, the restriction of $g'$ to $A$.
Then $\lambda_\omega(g)$ is the unique translation of $A$ such that $gv=\lambda_\omega(g)v$, for all $v\in S$. As the notation suggests, $\lambda_\omega(g)$ depends on $g$ and $\omega$, but not on $S$.

It is easy to check that the mapping $\lambda_\omega:g\mapsto\lambda_\omega(g)$ is a homomorphism from $B$ onto the group $T_0$ of type preserving translations on $A$.  Since $T_0 \cong \ZZ^n$ is an amenable group, it will follow that $B$ is amenable if $\ker\lambda_\omega$ can be shown to be amenable.

For each vertex $v$ of $[s,\omega)$, let $B_v=\{g\in B:gv=v\}$.  Then 
\[
\ker\lambda_\omega=\bigcup_{v\in[s,\omega)}B_v\, .
\]
Each of the groups $B_v$ is compact, being a closed subgroup of a parahoric subgroup.  
The group $\ker\lambda_\omega$ may thus be expressed as the inductive limit of the family of compact groups
$\{B_v:v\in[s,\omega)\}$, directed by inclusion.  Therefore $\ker\lambda_\omega$ is amenable.
\end{proof}

\begin{remark}
If $G$ is the group of $k$-rational points of an absolutely almost simple, simply connected linear algebraic $k$-group, this result is well known. For then the Borel subgroup $B$ is solvable, hence amenable.
\end{remark}

The amenability of the $\Gamma$-space $\Omega$ has the consequence that the Baum-Connes conjecture, with coefficients in $C(\Omega)$ has been verified \cite[Th\'eor\`eme 0.1]{Tu}.  Consequently $K_*(\cA_\Gamma)$ can be calculated by means of the Kasparov-Skandalis spectral sequence \cite[5.6, 5.7]{ks}.  This has initial terms
\begin{eqnarray*}
E^2_{p,q}
    &=& H_p(\Gamma,K_q(C(\Omega)))\\
    &=&\begin{cases} H_p(\Gamma,C(\Omega,\ZZ))\,,&\text{if $0\leq p\leq n$ and $q$ is even,}\\
                0\,,                &\text{otherwise.}
      \end{cases}
\end{eqnarray*}
Note that $H_p=0$ for $p>n$, since $\Gamma$ has homological dimension $\leq n$ .  Moreover $K_1(C(\Omega))=0$, since $\Omega$ is totally disconnected.

Suppose now that $\Delta$ has dimension $n=2$. 
Some of the nonzero terms in the first quadrant are shown in (\ref{spec2}).
\begin{equation}\label{spec2}
\begin{array}{lllllll}
 \cdot & \cdot & \cdot & \cdot & \cdot & \cdots \\ 
 E_{04}^2 & E_{14}^2 & E_{24}^2 & 0 & 0 & \cdots \\
  0 & 0 & 0 & 0 & 0 & \cdots \\
 E_{02}^2 & E_{12}^2 & E_{22}^2 & 0 & 0 & \cdots \\
  0 & 0 & 0 & 0 & 0 & \cdots \\
 E_{00}^2 & E_{10}^2 & E_{20}^2 & 0 & 0 & \cdots \\
\end{array}
\end{equation}
Recall that for $r\ge 2$ there are differentials $d_{p,q}^r : E^r_{p,q}\to E^r_{p-r,q+r-1}$,
and $E^{r+1}_{p,q}$ is the homology of $E_*^r$ at the position of $E^r_{p,q}$.
Since the differentials $d^2$ go up one row, it is clear that $d^2=0$ and $E^3_{p,q}=E^2_{p,q}$. Since the differentials $d^3$ go three units to the left,
$d^3=0$ and $E^4_{p,q}=E^3_{p,q}$. Continuing in this way we see that  
$E^\infty_{p,q}=E^2_{p,q}$.
Therefore the spectral sequence degenerates with $E^\infty_{p,q}=E^2_{p,q}$. Convergence of the spectral sequence to $K_*(\cA_\Gamma)$ means that
$$K_1(\cA_\Gamma)=H_1(\Gamma, C(\Omega,\ZZ))$$
and that there is a short exact sequence
$$0\longrightarrow H_0(\Gamma,C(\Omega,\ZZ))\longrightarrow K_0(\cA_\Gamma)
            \longrightarrow H_2(\Gamma, C(\Omega,\ZZ))\longrightarrow 0\,.$$
In particular, $\Omega_\Gamma=H_0(\Gamma,C(\Omega,\ZZ))$ is isomorphic to a subgroup of  $K_0(\cA_\Gamma)$.

\vskip 1cm

\section{$\widetilde A_2$ buildings and reduced group $C^*$-algebras}\label{a2tilde}

The reduced group $C^*$-algebra of a group $\G$ is the completion $C^*_r(\G)$ of the complex group algebra of $\G$ in the regular representation as operators on $\ell^2(\G)$. 
Let $\Gamma$ be a discrete torsion free group acting properly on the affine building $\Delta$, satisfying the hypotheses of Theorem \ref{main}.
By Proposition \ref{amenable}, $\Gamma$ acts amenably on the compact space $\Omega$. It follows
that the Baum-Connes assembly map is injective \cite{hig} and so the Novikov conjecture is true.
(This also follows from \cite{ks}.)
Therefore the class $[\1]$ in $K_0(C^*_r(\Gamma))$ does not have finite order.

Since $C_r^*(\G)$ embeds in $\cA_\G$, there is a natural homomorphism 
$$K_*(C_r^*(\G)) \to K_*(\cA_\Gamma).$$ 
This homomorphism is not injective, by Theorem \ref{main}, since $[\1]$ does have finite order in
$K_0(\cA_\Gamma)$. 
It is therefore worth comparing the $K$-theories of these two algebras. If the building is type $\widetilde A_2$, everything can be calculated explicitly.

The computation required is a corollary of \cite{la}, which states that the Baum-Connes conjecture
holds for any discrete group $\G$ satisfying the following properties.

\begin{itemize}
\item[(1)]
$\G$ acts continuously, isometrically and properly with compact quotient on a uniformly locally finite affine building or on a complete riemannian manifold of nonpositive curvature;
\item[(2)]
$\G$ has property (RD) of Jolissaint.
\end{itemize}
For a group $\G$ satisfying these conditions, $K_*(C_r^*(\G))$ is isomorphic to the geometric group 
$K_*^{\G}(\D) = KK_*^{\G}(C_0(\D),\CC)$. (The notation is consistent with \cite{bch}, because $\D$ is $\G$-compact.) This provides a way of calculating the groups $K_*(C_r^*(\G))$.

Assume therefore that all the conditions of Theorem \ref{main} hold, together with the condition that 
$\Delta$ has type $\widetilde A_2$. This is the case, for example, if $\Gamma$ is a torsion free lattice in $G=\SL_3(k)$.

Condition (1) is clearly satisfied and condition (2) is also satisfied by the main result of \cite{rrs}.
The finite cell complex $B\G= \D/\G$ is a $K(\G,1)$ space \cite[I.4]{bro}, so the group homology $H_*(\G,\ZZ)$ is isomorphic to the usual simplicial homology $H_*(B\G)$ \cite[Proposition II.4.1]{bro}.
Thus $H_0(\G,\ZZ)=\ZZ$ and $H_1(\G,\ZZ)=\G_{ab}$, the abelianization of $\G$.
Moreover, since $B\G$ is $2$-dimensional, the group $\G$ has homological dimension at most $2$ \cite[VIII.2 Proposition (2.2) and VIII.6 Exercise 6]{bro}. It follows that  $H_2(\G,\ZZ)$ is free abelian and $H_p(\G,\ZZ)=0$ for $p>2$. Since $\G$ satisfies the Baum-Connes conjecture, $K_*(C_r^*(\G))$ coincides with its ``$\gamma$-part'' \cite[Definition-corollary 3.12]{ka}. Therefore $K_*(C_r^*(\G))$ may be computed as the limit of a spectral sequence $E^r_{p,q}$ \cite [Theorem 5.6 and Remark 5.7(a)]{ks}.
Since $\G$ is torsion free, $\G$ acts freely on $\D$. 
According to \cite [Remarks 5.7(b)]{ks} the initial terms of the spectral sequence are
\begin{equation}\label{Epq}
E^2_{p,q}= H_p(\G, K_q(\CC))= \begin{cases}
H_p(\G,\ZZ) & \text{if $p \in \{0,1,2\}$ and $q$ is even},\\
0& \text{otherwise.}
\end{cases}
\end{equation}
The nonzero terms in the first quadrant are shown in (\ref{spec2}).
Exactly as for (\ref{spec2}), the spectral sequence degenerates with $E^\infty_{p,q}=E^2_{p,q}$.
Convergence of the spectral sequence to $K_*(C_r^*(\G))$ means that
\begin{equation}\label{K1cstar}
K_1(C_r^*(\G))= H_1(\G,\ZZ)
\end{equation}
and that there is a short exact sequence
\begin{equation}\label{exactcstar}
0\longrightarrow H_0(\G,\ZZ)\longrightarrow K_0(C_r^*(\G))\longrightarrow H_2(\G,\ZZ)\longrightarrow 0.
\end{equation}
Since $H_2(\G,\ZZ)$ is free abelian, (\ref{exactcstar}) splits and we have
\begin{equation}\label{kgp}
\begin{split}
K_0(C_r^*(\G)) &= H_0(\G,\ZZ)\oplus H_2(\G,\ZZ),\\
K_1(C_r^*(\G)) &= H_1(\G,\ZZ).
\end{split}
\end{equation}
Now $H_1(\G,\ZZ)$ is a finite group, because $\G$ has Kazhdan's property (T) \cite[Corollary 1]{bs}. 
It follows that $K_0(C_r^*(\G)) = \ZZ^{\chi(\G)}$, where $\chi(\G)$ is the Euler-Poincar\'e characteristic of $\G$. This proves

\begin{theorem}\label{main2}Let $\G$ be a torsion free cocompact lattice in $G$, where
$(G,{\mathfrak I},N,S)$ is an affine topological Tits system of type  $\widetilde A_2$. 
Then
\begin{equation}\label{main1}
K_0(C_r^*(\G)) = \ZZ^{\chi(\G)} 
\quad \text{and}\quad
K_1(C_r^*(\G)) = \G_{ab}.
\end{equation}
\end{theorem}

\medskip

The value of $\chi(\G)$ is easily calculated \cite[p. 150, Th\'eor\`eme 7]{ser0},\cite[Section 4]{R}. It is 
\begin{equation}\label{ep}
\chi(\G)=(q-1)(q^2-1) \cdot n_\fs(\Gamma),
\end{equation}
where $q$ is the order of the building $\D$ and $n_\fs(\Gamma)$ is the number of $\G$-orbits of vertices of 
type $\fs$, where $\fs\in I$ is fixed.

\bigskip

In \cite{cmsz} a detailed study was undertaken of groups of type rotating automorphisms of $\widetilde A_2$ buildings, subject to the condition that the group action is free and transitive
on the vertex set of the building. For $\widetilde A_2$ buildings of orders $q=2,3$, the authors of that article give a complete enumeration of the possible groups with this property. These groups are
called $\widetilde A_2$ groups. Some, but not all, of the $\widetilde A_2$ groups are cocompact lattices in $\pgl_3(k)$ for some 
local field $k$ with residue field of order $q$. It is an empirical fact 
that either $k=\QQ_p$ or $k= \FF _q((X))$ in all the examples constructed so far.

For each $\widetilde A_2$ group $\widetilde \Gamma < \pgl_3(k)$, consider the unique type preserving subgroup $\Gamma < \widetilde \Gamma$ of index 3.
Each such $\Gamma$ is torsion free and acts freely and transitively on the 
set of vertices of a fixed type $\fs$. That is $n_\fs=1$. Therefore
\[
\chi(\G)=(q-1)(q^2-1)=1+ \rank H_2(\G,\ZZ)\, .
\] 

\begin{remark}\label{zerot}
There are eight such groups $\Gamma$ if $q=2$, and twenty-four if $q=3$.
Using the results of \cite{RS} and the MAGMA computer algebra package, one can compute
$K_0(\cA_\Gamma)$. 
One checks that in all these examples,
\begin{equation*}
\rank  K_0(\cA_\Gamma) = 2\cdot\rank H_2(\G,\ZZ) = 
\begin{cases}
    4&  \text{if $q=2$},\\
    30 &  \text{if $q=3$.}
\end{cases}
\end{equation*}
Furthermore, the class of $[\1]$ in the $K_0(\cA_\Gamma)$ has order $q-1$.
Note that for $q=2$ this means that $[\1]=0$.
\end{remark}

These values also appear to be true for higher values of $q$. In particular, they have been verified for a number of groups with $q=4,5,7$.  Here is an example with $q=4$.

\begin{example}\label{ex1.1}
Consider the Regular $\widetilde A_2$ group $\Gamma_r$, with $q=4$. This is a torsion free cocompact
subgroup of $\pgl_3(\KK)$, where $\KK$ is the Laurent series field 
${\boldmath F}_4((X))$ with coefficients in the field ${\boldmath F}_4$ with four elements. It is described in \cite[Part I, Section 4]{cmsz},
and its embedding in $\pgl_3({\boldmath F}_4((X)))$ is essentially unique,
by the Strong Rigidity Theorem of Margulis.
The group $\Gamma_r$ is torsion free and has 21 generators
$x_i, 0\le i \le 20$, and relations (written modulo 21):
$$
\begin{cases}
x_jx_{j+7}x_{j+14}=x_jx_{j+14}x_{j+7}=1 \quad & 0\le j\le 6,\\
x_jx_{j+3}x_{j-6}=1 \quad & 0\le j\le 20.
\end{cases}
$$
Let $\Gamma<\psl_3(\KK)$ be the type preserving index three subgroup of $\Gamma_r$.
The group $\Gamma$ has generators $x_jx_0^{-1}$, $1\le j\le 20$.
Using the results of \cite{RS} one obtains 
$$K_0(\cA_\Gamma)=\ZZ^{88}\oplus (\ZZ/2\ZZ)^{12}\oplus (\ZZ/3\ZZ)^4\oplus(\ZZ/7\ZZ)^4\oplus(\ZZ/9\ZZ),$$
and the class of $[\1]$ in $K_0(\cA_\Gamma)$ is $3+\ZZ/9\ZZ$, which has order $q-1=3$. It also follows from
\cite[Theorem 2.1]{RS} that $K_0(\cA_\Gamma)=K_1(\cA_\Gamma)$.

According to Theorem \ref{main2},
\begin{equation*}
K_0(C_r^*(\G)) = \ZZ^{45} =\ZZ^{44}\oplus\langle [\1]\rangle 
\quad \text{and}\quad
K_1(C_r^*(\G)) = (\ZZ/2\ZZ)^6\oplus (\ZZ/3\ZZ).
\end{equation*}
The second equality was obtained using the MAGMA computer algebra package.
This, and similar, examples suggest that the {\it only} reason for failure of injectivity of the
natural homomorphism 
$$K_0(C_r^*(\G)) \to K_0(\cA_\Gamma)$$
is the fact that $[\1]$ has finite order in $K_0(\cA_\Gamma)$.
\end{example}

\begin{example}\label{ex1.2}
For completeness, here are the results of the computations for one of the groups with $q=3$.
The Regular group ${1.1}$ of \cite{cmsz}, with $q=3$, has 13 generators
$x_i$, $0\le i \le 12$, and relations (written modulo 13):
$$
\begin{cases}
x_j^3=1 \quad & 0\le j\le 13,\\
x_jx_{j+8}x_{j+6}=1 \quad & 0\le j\le 13.
\end{cases}
$$
Let $\Gamma$ be the type preserving index three subgroup.
The group $\Gamma$ has generators $x_jx_0^{-1}$, $1\le j\le 12$.
Note that the group $1.1$ has torsion, but its type preserving subgroup $\Gamma$ is torsion free.
One obtains 
$$K_0(\cA_\Gamma)=\ZZ^{30}\oplus (\ZZ/2\ZZ)\oplus (\ZZ/3\ZZ)^6\oplus(\ZZ/13\ZZ)^4,$$
and the class of $[\1]$ in $K_0(\cA_\Gamma)$ is $1+\ZZ/2\ZZ$, which has order $q-1=2$. It also follows 
from Theorem \ref{main2} that
\begin{equation*}
K_0(C_r^*(\G)) = \ZZ^{16} 
\quad \text{and}\quad
K_1(C_r^*(\G)) = (\ZZ/3\ZZ)^3\oplus (\ZZ/13\ZZ).
\end{equation*}
\end{example}

\end{document}